\documentclass[a4paper, 12pt] {article}
\usepackage[cp1251]{inputenc}
\usepackage[english]{babel}
\usepackage{amsfonts}
\usepackage{mathtext}
\usepackage{xcolor}
\usepackage{cite}

\begin{document}
\thispagestyle{empty}
\begin{center}
{\bf ON DAMPING A CONTROL SYSTEM WITH GLOBAL AFTEREFFECT ON QUANTUM GRAPHS. STOCHASTIC INTERPRETATION }
\end{center}

\begin{center}
{\large\bf Sergey Buterin}
%\footnote{Department of Mathematics, Saratov University, Saratov, Russia, {\it Email: buterinsa@sgu.ru}}
\end{center}

{\bf Abstract.} Quantum graphs model processes in complex systems represented as spatial networks in various
fields of natural science and technology. An example is the oscillations of elastic string networks, the
nodes of which, besides the continuity conditions, also obey the Kirchhoff conditions, expressing the balance
of tensions. In this paper, we propose a new look at quantum graphs as {\it temporal} networks, which means
that the variable parametrizing the edges of a graph is interpreted as time, while each internal vertex is a
branching point giving several different scenarios for the further trajectory of a process. Then
Kirchhoff-type conditions may also arise. Namely, they will be satisfied by such a trajectory of the process
that is optimal with account of all the scenarios simultaneously. By employing the recent concept of global
delay, we extend the problem of damping a first-order control system with aftereffect, considered earlier
only on an interval, to an arbitrary tree graph. The first means that the delay, imposed starting from the
initial moment of time, associated with the root of the tree, propagates through all internal vertices.
Bringing the system into the equilibrium and minimizing the energy functional with account of the anticipated
probability of each scenario, we come to a variational problem. Then, we establish its equivalence to a
self-adjoint boundary value problem on the tree for some second-order equations involving both the global
delay and the global advance. The unique solvability of both problems is proved. We also illustrate that the
interval case when the coefficients of the equation are discrete stochastic processes in discrete time can be
viewed as the extension to a tree.


\\

{\bf 1. Introduction}
\\

Purely differential (local) operators on graphs, often called {\it quantum graphs}, model various processes
in {\it spatial} networks appearing in organic chemistry, mesoscopic physics, quantum mechanics,
nanotechnology, hydrodynamics, waveguide theory, and many other fields \cite{Mont, Nic, vB, LangLeug, Kuch,
BCFK, BerkKuch, Pok, Kuz-17, Boris-22}. Among the simplest applications, one can mention elastic string
networks, whose nodes obey continuity conditions expressing the firm junction of the strings, as well as
Kirchhoff conditions expressing, in turn, the balance of tensions. Such matching conditions at the internal
vertices and their generalizations are also relevant to other models dealing with spatial quantum graphs.

In the present paper, we propose a new look at quantum graphs as {\it temporal} networks. This means that the
variable parametrizing the edges will be associated with time, while at each internal vertex, there will
appear several different scenarios for the further course of a process in accordance with the number of edges
emanating from that vertex. Such settings immediately imply only the continuity matching conditions at the
internal vertices. However, the interesting effect observed below is that Kirchhoff-type conditions remain
relevant but to such a trajectory of the process that is optimal with account of all the scenarios
simultaneously.

Specifically, we extend to graphs the problem of damping a control system with aftereffect described by an
equation of retarded type, which was posed and studied on an interval by Krasovskii in~\cite{Kras-68}. Later
on, Skubachevskii \cite{Skub-94} considered a generalization of this problem to the case in which the
equation contains the delay also in the dominant terms, i.e. belongs to neutral type:
\begin{equation}\label{1.1}
y'(t)+ay'(t-\tau)+by(t)+cy(t-\tau)=u(t), \quad t>0,
\end{equation}
where $a,\,b,\,c\in{\mathbb R}$ and $\tau>0$ are fixed, while $u(t)$ is a real-valued square-integrable
control function. Recall that equation (\ref{1.1}) is said to be of {\it neutral type} if $a\ne0.$ Otherwise,
it belongs to {\it retarded type} (if $c\ne0).$ A previous history of the system is determined by the
condition
\begin{equation}\label{1.2}
y(t)=\varphi(t), \quad t\in[-\tau,0],
\end{equation}
where $\varphi(t)$ is a known real-valued function in $W_2^1[-\tau,0],$ i.e. $\varphi(t)$ is absolutely
continuous on $[-\tau,0]$ and $\varphi'(t)\in L_2(-\tau,0).$

One can show by steps that the system (\ref{1.1}), (\ref{1.2}) has a unique solution $y(t)\in W_2^1[0,T]$ for
any $T>0.$ Indeed, relations (\ref{1.1}) and (\ref{1.2}) give the Cauchy problem
$$
y'(t)+by(t)=u(t)-a\varphi'(t-\tau)-c\varphi(t-\tau), \quad 0<t<\tau, \quad y(0)=\varphi(0),
$$
having a unique solution $y(t)\in W_2^1[0,\tau].$ By continuing the process, the solution $y(t)$ can be
successively obtained on the subsequent intervals $(j\tau,(j+1)\tau]$ for all $j\in{\mathbb N}.$

Fix some $T>\tau.$ The control problem can be formulated as follows. Find $u(t)\in L_2(0,T)$ that would bring
the system (\ref{1.1}), (\ref{1.2}) into the equilibrium $y(t)=0$ for $t\ge T.$

For achieving this, it is sufficient to find $u(t)$ leading to the state
\begin{equation}\label{1.3}
y(t)=0, \quad t\in[T-\tau,T],
\end{equation}
and put $u(t)=0$ for $t>T.$ Since such $u(t)\in L_2(0,T)$ is not unique whenever $T>\tau,$ it is reasonable
to look for it trying to minimize the efforts $\|u\|_{L_2(0,T)}.$

Thus, there arises the variational problem for the energy functional
\begin{equation}\label{1.4}
J(y)=\int_0^T \left(y'(t)+ay'(t-\tau)+by(t)+cy(t-\tau)\right)^2\,dt \to\min
\end{equation}
under the conditions (\ref{1.2}) and (\ref{1.3}).

Solution of the problem (\ref{1.2})--(\ref{1.4}) was obtained in \cite{Skub-94} and also given in the
monograph~\cite{Skub-97}. In particular, it has been reduced to a self-adjoint boundary value problem for a
second-order functional-differential equation, which involves both negative and positive shifts $\pm\tau$ of
the argument as long as $T>2\tau.$ Specifically, the following theorem holds.

\medskip
{\bf Theorem 1. \cite{Skub-94} }{\it A function $y(t)\in W_2^1[-\tau,T]$ is a solution of the variational
problem (\ref{1.2})--(\ref{1.4}) if and only if $y(t)$ solves the boundary value problem for the equation
$$
\left((1+a^2)y'(t) +ay'(t-\tau)+ay'(t+\tau)\right)' +(c-ab)(y'(t-\tau)-y'(t+\tau))
$$
\begin{equation}\label{1.4-1}
=(b^2+c^2)y(t) +bc(y(t-\tau)+y(t+\tau)), \quad 0<t<T-\tau,
\end{equation}
under the conditions (\ref{1.2}) and (\ref{1.3}).}

\medskip
The boundary value problem (\ref{1.2}), (\ref{1.3}), (\ref{1.4-1}) may have no solution in $W_2^2[0,T-\tau]$
if $a\ne0$ even when $\varphi(t)\in W_2^2[-\tau,0].$ For this reason, the solution is understood in the
generalized sense, i.e.
$$
(1+a^2)y'(t) +ay'(t-\tau)+ay'(t+\tau)\in W_2^1[0,T-\tau].
$$
The following theorem gives the existence and uniqueness of such a generalized solution, which also means the
unique solvability of the variational problem (\ref{1.2})--(\ref{1.4}).

\medskip
{\bf Theorem 2. \cite{Skub-94} }{\it For any function $\varphi(t)\in W_2^1[-\tau,0],$ there exists a unique generalized solution $y(t)\in
W_2^1[-\tau,T]$ of the problem (\ref{1.2}), (\ref{1.3}), (\ref{1.4-1}). Moreover, $y(t)$ obeys the estimate
$$
\|y\|_{W_2^1[-\tau,T]}\le c\|\varphi\|_{W_2^1[-\tau,0]},
$$
where $c$ does not depend on $\varphi(t).$}

\medskip
In \cite{Skub-94}, also necessary and sufficient conditions in terms of $\varphi(t)\in W_2^k[-\tau,0]$ were
obtained for the generalized solution $y(t)$ to belong to the space $W_2^k[-\tau,T]$ for $k\ge2.$

For generalizations as well as analogs of this problem on an interval, see more recent papers \cite{Ross-17,
AkhdSkub, IvanSkub} and the references therein. Various optimal control problems and variational problems
with deviating argument on an interval can be found also in \cite{ Hal-68, Kent-71, BaKe-72, BaJa-73,
KharTad-78, Col-88, BravDrakh-95} and other works.

The systematic study of equations with delay began in the second third of the last century due to the growing
number of applications \cite{Mysh}. In particular, both scalar and vectorial equations of the form
(\ref{1.1}) have been extensively studied (see also \cite{BellCook}). There are a number of physical
situations which motivate the control problems as above (see, e.g., \cite{BaKe-72}). For example, let $y(t)$
represent some error which one wishes to be driven to zero and held there if possible. Then if the error is
described by (\ref{1.1}), then it is obvious that the desired terminal condition is (\ref{1.3}). Another
interesting example given in \cite{BaKe-72} deals with a boundary control of linear hyperbolic partial
differential equations, which can be transformed to problems involving control of functional-differential
equations.

\medskip
An extension of the above control problem to graphs requires an appropriate definition of the
functional-differential equation with delay on them. However, various functional-differential as well as
other classes of nonlocal operators on graphs until recently were considered only in the {\it locally}
nonlocal case, when the corresponding nonlocal equation on each edge can be solved independently of the
equations on the other edges \cite{Nizh-12, Bon18-1, HuBondShYan19, Hu20, WangYang-21, Bon22, WangYang-22}.
This limitation always left unclear how the nonlocalities could coexist with the internal vertices of the
graph and, in particular, how one could describe a process with global aftereffect on the entire graph.

This gap was addressed in \cite{But23}, where a concept of functional-differential operators on graphs with
global delay was suggested. This means that the delay, being measured in the direction of a specific boundary
vertex, called the root, propagates through all internal vertices of the graph to the subsequent edges. In
particular, for each internal vertex of a directed tree graph, a solution of the equation defined on its
incoming edge serves as an initial function for the delayed equations defined on the outgoing edges. This
idea naturally comes from the equation on a plain interval, which can be formally divided into two parts
(edges) by any internal point. Then the solution on the first part automatically becomes an initial function
for the same equation on the second one.

The global delay has become an alternative to the locally nonlocal settings introduced by Wang and Yang in
\cite{WangYang-21}, where the equation on each edge possessed its own delay parameter and could be considered
independently. The papers \cite{But23, But23-2, WYBD24} deal with some inverse spectral problems in the
highlighted globally nonlocal case. Meanwhile, this concept also appears natural for the extension of
equation (\ref{1.1}) along with the corresponding control problem to graphs.

\medskip
For simplicity, we focus here on the case $a=0,$ returning thus to equation (\ref{1.1}) of retarded~type.
However, we pay special attention to analyzing the extension of the original control problem to a tree. In
particular, we show that in the structure of the corresponding energy functional, it is natural to take into
account the probabilities of all the scenarios arising at its internal vertices.

The case $a\ne0$ on graphs as well as a development of the concept of generalized solution from
\cite{Skub-94} are addressed in \cite{But24} within the more general situation involving a control system of
an arbitrary-order and neutral type with non-smooth complex coefficients.

Following in general the strategy for an interval, we establish equivalence of the corresponding variational
problem to a certain self-adjoint boundary value problem on the tree for second-order functional-differential
equations with bidirectional shifts of the argument. Then, the unique solvability of both problems is proved.
An estimate for the solution via the initial function determining the prehistory of the process on the entire
tree is also obtained. However, the graph case essentially complicates the study. In particular, there
appears new Lemma~2, which is crucial for deriving the boundary value problem on a tree. That problem, in
turn, acquires complicated globally nonlocal structure including Kirchhoff-type conditions. To the best of
our knowledge, the latter previously arose only in {\it spatial} networks, and their emergence for a time
variable becomes a new quality.

In the next section, we illustrate the obtained results on a star-shaped graph and discuss some possible
interpretations of this problem on a graph. An interesting interpretation, which will be further expanded in
the last section to arbitrary trees, brings a stochastic nature into the original equation~(\ref{1.1}).
Specifically, the system on a tree arises when one replaces the constant coefficients in~(\ref{1.1}) with
discrete-time stochastic processes having a finite number of states in ${\mathbb R}.$ Then a decomposition of
those processes into all their possible realizations leads to a system of equations with constant
coefficients but defined on a finite tree (see Section~6 for details).

\medskip
The paper is organized as follows. In Section~2, we consider the case of a star-shaped graph. In Section~3,
we introduce the control system with aftereffect on an arbitrary tree and formulate the corresponding
variational problem. In Section~4, we establish its equivalence to a boundary value problem on the tree. The
unique solvability of both problems will be proved in Section~5. In the last section, we discuss the
stochastic interpretation of the control problem on a tree.
\\

{\bf 2. Star-shaped graph}
\\

Let up to the time point $t=T_1$ associated with the internal vertex $v_1$ of the graph $\Gamma_m$ in Fig.~1,
our control system with delay $\tau<T_1$ be described by the equation
\begin{equation}\label{1.5}
\ell_1y(t):=y_1'(t)+b_1y_1(t)+c_1y_1(t-\tau)=u_1(t), \quad 0<t<T_1,
\end{equation}
where $y_1(t)$ is defined on the edge $e_1=[v_0,v_1]$ of $\Gamma_m$ and has the prehistory
\begin{equation}\label{1.9}
y_1(t)=\varphi(t), \quad t\in[-\tau,0].
\end{equation}
At the vertex~$v_1,$ this system branches into $m-1$ independent processes described by the equations
\begin{equation}\label{1.6}
\ell_jy(t):=y_j'(t)+b_jy_j(t)+c_jy_j(t-\tau)=u_j(t), \quad t>0, \quad j=\overline{2,m},
\end{equation}
but having a common history determined by equation (\ref{1.5}) along with the conditions (\ref{1.9}) and
\begin{equation}\label{1.7}
y_j(t)=y_1(t+T_1), \quad t\in(-\tau,0), \quad j=\overline{2,m}.
\end{equation}
Besides (\ref{1.7}), it is natural to impose continuity conditions at~$v_1:$
\begin{equation}\label{1.8}
y_j(0)=y_1(T_1), \quad j=\overline{2,m},
\end{equation}
agreeing with (\ref{1.7}) as $t\to0^-$ (Remark~1). Here and below, $j=\overline{j_1,j_2}$ means
$j=j_1,j_1+1,\ldots, j_2.$

 As in the preceding section, we assume that $b_j,c_j\in {\mathbb
R},\;j=\overline{1,m}.$

\begin{center}
\unitlength=0.75mm
\begin{picture}(80,107)
\multiput(-10,56)(-1,0){14}{\circle*{0.7}} \put(-23.4,56){\circle*{1}}                                   \put(-29,58){\small $-\tau$}

\put(-10,56){\line(1,0){50}}    \put(11,60){\small $e_1$} \put(-10,56){\circle*{1}}       \put(-13,50){\small
$v_0$}    \put(-12.7,58){\small $0$}
                                                              \put(38,63){\small $0$}
                                                              \put(42.5,60){\small $0$}
                                                              \put(42,49){\small $0$}
\put(26.6,56){\circle*{1}}                                    \put(11.5,50){\small $T_1-\tau$} \put(40,56){\circle*{1}}
\put(33.5,50){\small $v_1$}     \put(33,58){\small $T_1$}

\put(40,56){\vector(1,3){15}}     \put(55,102){\small $e_2$} \put(47,77){\circle*{1}} \put(47,72){\small
$T_2-\tau$} \put(51.3,89.9){\circle*{1}} \put(43,91){\small $T_2$}

\put(40,56){\vector(2,1){61}}     \put(101,84){\small $e_3$} \put(78,75){\circle*{1}}
\put(78.5,70.5){\small $T_3-\tau$} \put(90,81){\circle*{1}}                                      \put(86,84){\small $T_3$}

\multiput(40,56)(2,-1){22}{\circle*{0.7}}

\put(40,56){\vector(1,-3){13}}     \put(52.5,14.5){\small $e_m$} \put(45,41){\circle*{1}}
\put(25.5,37){\small $T_m-\tau$} \put(49.3,28.1){\circle*{1}}
\put(50,27){\small $T_m$}

\put(1,0){\small Fig. 1. A star-shaped graph $\Gamma_m$}
\end{picture}
\end{center}

For $j=\overline{2,m},$ the $j$th equation in (\ref{1.6}) is defined on the edge $e_j$ of $\Gamma_m,$ which
is originally an infinite ray emanating from the internal vertex $v_1.$ Conditions (\ref{1.7}) mean that the
delay propagates through $v_1.$

The problem (\ref{1.5})--(\ref{1.8}) has a unique solution $y_j(t)\in W_2^1[0,T_j],\;j=\overline{1,m},$ for
any fixed $T_j>0,$ $j=\overline{2,m},$ whenever $u_j(t)\in L_2(0,T_j)$ for $j=\overline{1,m}.$ Indeed, on the
edge $e_1,$ it can be solved as the Cauchy problem (\ref{1.1}), (\ref{1.2}). Then, the obtained $y_1(t)$
gives the common initial function (\ref{1.7}) and the common initial condition (\ref{1.8}) for all equations
in (\ref{1.6}), which can be solved similarly.

\medskip
{\bf Example 1.} Let $m=2,$ $b:=b_1=b_2,$ $c:=c_1=c_2$ and
$$
y(t):=\left\{\begin{array}{cc}y_1(t), & 0\le t\le T_1,\\[3mm]
y_{2}(t-T_1), & t>T_1,
\end{array}\right.
 \quad
u(t):=\left\{\begin{array}{cc}u_1(t), & 0< t< T_1,\\[3mm]
u_{2}(t-T_1), & t>T_1.
\end{array}\right.
$$
Then the Cauchy problem (\ref{1.5})--(\ref{1.8}) takes the form (\ref{1.1}), (\ref{1.2}) with $a=0.$

\medskip
For definiteness, let $T_j>\tau$ also for all $j=\overline{2,m}.$ Analogously to the preceding section, we
assume that one needs to bring the system (\ref{1.5})--(\ref{1.8}) into the equilibrium state
\begin{equation}\label{1.10}
y_j(t)=0, \quad t\in[T_j-\tau,T_j], \quad j=\overline{2,m},
\end{equation}
by choosing suitable controls $u_j(t)\in L_2(0,T_j),\;j=\overline{1,m}.$ Then, letting $u_j(t)=0$ for $t>T_j$
and $j=\overline{2,m}$ will guarantee $y_j(t)=0$ for such $t$ and $j.$ In other words, the system will be
damped on each outgoing edge. Since such $u_j(t)$ are not unique, it is natural to try reducing the efforts
$\|u_j\|_{L_2(0,T_j)}$ as much as possible. One can also regulate the participation of each
$\|u_j\|_{L_2(0,T_j)}^2$ in the corresponding energy functional by choosing a certain positive weight
$\alpha_j.$

Specifically, we consider the variational problem
\begin{equation}\label{1.11}
\sum_{j=1}^m \alpha_j\int_0^{T_j}(\ell_jy(t))^2\,dt\to\min
\end{equation}
under the conditions (\ref{1.9}) and (\ref{1.7})--(\ref{1.10}), where  $\alpha_j>0,\;j=\overline{1,m},$ are fixed.

Let us discuss some possible interpretations of the control system on $\Gamma_m.$

\medskip
{\bf (i)} The first possibility is that the process really splits into $m-1$ independent processes at the
time point $t=T_1.$ This means that the initial process by some means gives rise to such independent
processes described by their own equations differing in the coefficients $b_j$ and $c_j.$ It is required, in
turn, that each of the arisen processes should be  damped starting from the corresponding time point $t=T_j.$
Since all the processes really exist and if there is no reason to emphasize certain ones from the others, it
may be advisable to take equal $\alpha_j$ in (\ref{1.11}).

\medskip
{\bf (ii)} Another interesting interpretation occurs when not all of the ''arisen'' processes will really
unfold. For example, all of them are just possible scenarios for one and the same process after the time
point $t=T_1.$ In other words, at $t=T_1,$ there appear $m-1$ different scenarios of the further process flow
determined, in turn, by different pairs of the coefficients $b_j$ and $c_j$ in equations (\ref{1.6}). Before
$t=T_1,$ there is no information about which scenario will~be~really~fulfilled.  However, the system should
be surely damped at each possible outcome.

\medskip
Under the second interpretation, it is reasonable to choose the weights $\alpha_j$ in (\ref{1.11}) according to the anticipated probability
of the corresponding scenario. For example, one can take
\begin{equation}\label{1.11-1}
\alpha_1=1, \quad \alpha_j=\frac1{m-1}, \quad j=\overline{2,m},
\end{equation}
bearing in mind that the probability of the very first scenario (when $t<T_1)$ always equals~$1,$ and also
assuming all further scenarios starting from $t=T_1$ to be equiprobable.

Before continuing this discussion and giving some illustrative examples, let us first formulate our results on the solvability of the
variational problem (\ref{1.9}), (\ref{1.7})--(\ref{1.11}).

The following assertion is a particular case of Theorem~5 for an arbitrary tree.

\medskip
{\bf Theorem 3. }{\it Functions $y_1(t)\in W_2^1[-\tau,T_1],$ $y_j(t)\in W_2^1[0,T_j],\,j=\overline{2,m},$ form a solution of the variational
problem (\ref{1.9}), (\ref{1.7})--(\ref{1.11}) if and only if they possess the additional smoothness
$$
y_1(t)\in W_2^2[0,T_1], \quad y_j(t)\in W_2^2[0,T_j-\tau], \quad j=\overline{2,m},
$$
and solve the boundary value problem (which we denote by ${\mathcal B})$ consisting of the equations
$$
\alpha_1(\ell_1 y)'(t)=\alpha_1b_1\ell_1 y(t)+\left\{\begin{array}{cc}
\alpha_1c_1\ell_1 y(t+\tau), & 0<t<T_1-\tau,\\[3mm]
\displaystyle\sum_{\nu=2}^m \alpha_\nu c_\nu\ell_\nu y(t+\tau-T_1), & T_1-\tau< t<T_1,
\end{array}\right.
$$
$$
(\ell_j y)'(t)=b_j\ell_j y(t)+ c_j\ell_j y_j(t+\tau), \quad 0<t<T_j-\tau, \quad j=\overline{2,m},
$$
along with the standing conditions (\ref{1.9}) and (\ref{1.7})--(\ref{1.10}) as well as the Kirchhoff-type
condition
\begin{equation}\label{1.12}
\alpha_1y_1'(T_1) +\Big(\alpha_1b_1 -\sum_{j=2}^m\alpha_jb_j\Big) y_1(T_1) +\Big(\alpha_1c_1
-\sum_{j=2}^m\alpha_jc_j\Big)y_1(T_1-\tau)=\sum_{j=2}^m \alpha_j y_j'(0),
\end{equation}
additionally emerging at the internal vertex $v_1.$}

\medskip
In particular, if all $\alpha_j$ are equal and the coefficients in (\ref{1.5}) and (\ref{1.6}) satisfy the relations
\begin{equation}\label{1.13}
b_1 =\sum_{j=2}^mb_j, \quad c_1 =\sum_{j=2}^mc_j,
\end{equation}
then (\ref{1.12}) becomes classical Kirchhoff's condition (see also Remark~2), which often arises in the
theory of spatial quantum graphs. For example, it expresses the balance of tensions in a system of connected
strings or Kirchhoff's law in electrical circuits. If the second expression in parentheses in (\ref{1.12})
does not vanish, then (\ref{1.12}) can be classified as a nonlocal Kirchhoff-type condition.

Relations (\ref{1.13}) generalize the constancy of the coefficients in (\ref{1.1}), while their absence would
correspond to the appearance of stepwise functions in (\ref{1.1}) instead of $b$ or~$c.$ In the latter case,
due to (\ref{1.12}), the solution of the resulting boundary value problem on an interval may lose the
smoothness at the discontinuities of such stepwise coefficients $b(t)$ and~$c(t)$ (see also Remark~5).

The next theorem, being a particular case of Theorem~6 for an arbitrary tree, also gives the unique
solvability of the variational problem (\ref{1.9}), (\ref{1.7})--(\ref{1.11}).

\medskip
{\bf Theorem 4. }{\it The boundary value problem ${\mathcal B}$ has a unique solution.
Moreover, this solution satisfies the estimate
$$
\|y_1\|_{W_2^1[0,T_1]}+\sum_{j=2}^m \|y_j\|_{W_2^1[0,T_j-\tau]} \le C\|\varphi\|_{W_2^1[-\tau,0]},
$$
where $C$ is independent of $\varphi(t).$}

\medskip
These theorems establish the existence and uniqueness of an optimal trajectory $[y_1,y_2\ldots,y_m]$ with
account of all possible scenarios $e_2,\ldots, e_m$ simultaneously. Substituting$[y_1,y_2\ldots,y_m]$ into
equations (\ref{1.5}) and (\ref{1.6}), one can obtain the corresponding optimal control
$[u_1,u_2,\ldots,u_m].$ Let us see how this control should be used in frames of the interpretation (ii)
outlined above.

Specifically, before the time point $t=T_1,$ i.e. in equation (\ref{1.5}), one should use~$u_1(t).$ Further,
at $t=T_1,$ it becomes known which scenario $j_0\in\{2,\ldots,m\}$ unfolds. Thus, we actually deal only with
the $j_0$th equation in (\ref{1.6}) and should use the corresponding control $u_{j_0}(t).$

Of course, the resulting composite control $[u_1,u_{j_0}]$ is, generally speaking, more ''expensive'' than
the one that could be obtained by optimizing the system solely along the timeline $[e_1,e_{j_0}]$ in case
when the fulfilment of the $j_0$th scenario were known a priori. However, $[u_1,u_{j_0}]$ appears to be the
best choice under the uncertainty when each scenario among $e_2,\ldots,e_m$ is possible.

\medskip
{\bf Example 2.} Let $b_j,\;c_j$ and $T_j$ be independent of $j\in\{2,\ldots,m\}.$ Then we have $m-1$ copies
of one and the same scenario starting from the time point $t=T_1.$ By the symmetry, the solution
$[y_1,y_2,\ldots, y_m]$ of the boundary value problem ${\cal B}$ under the assumption~(\ref{1.11-1}) contains
$m-1$ equal components: $y_2(t)\equiv\ldots\equiv y_m(t)$ (otherwise, this solution would not be unique).
Hence, the number of equations in the problem ${\cal B}$ can be reduced to just two, namely:
\begin{equation}\label{1.13-1}
(\ell_1 y)'(t)=b_1\ell_1 y(t)+\left\{\begin{array}{cc}
c_1\ell_1 y(t+\tau), & 0<t<T_1-\tau,\\[3mm]
\displaystyle c_{2}\ell_{2} y(t+\tau-T_1), & T_1-\tau< t<T_1,
\end{array}\right.
\end{equation}
and
\begin{equation}\label{1.13-2}
(\ell_{2} y)'(t)=b_{2}\ell_{2} y(t)+ c_{2}\ell_{2} y_{2}(t+\tau), \quad 0<t<T_{2}-\tau.
\end{equation}
Moreover, conditions (\ref{1.7})--(\ref{1.10}) will take the forms
\begin{equation}\label{1.13-3}
y_{2}(t)=y_1(t+T_1), \; t\in(-\tau,0);  \quad y_{2}(0)=y_1(T_1); \quad y_{2}(t)=0, \; t\in[T_{2}-\tau,T_{2}],
\end{equation}
respectively, while the Kirchhoff condition (\ref{1.12}) can be represented as
\begin{equation}\label{1.13-4}
y_1'(T_1) +(b_1 -b_{2}) y_1(T_1) +(c_1 -c_{2})y_1(T_1-\tau)= y_{2}'(0).
\end{equation}
In particular, when $b_1 =b_2=:b$ and $c_1 =c_2=:c,$ the problem (\ref{1.9}), (\ref{1.13-1})--(\ref{1.13-4})
coincides with the problem (\ref{1.2}), (\ref{1.3}), (\ref{1.4-1}) for $a=0,$ $T=T_1+T_{2}$ and $y(t)$
defined in Example~1.

\medskip
Thus, artificial reproducing copies of one and the same scenario starting from a certain point of the
interval and employing appropriate weights in the corresponding energy functional (\ref{1.11}) leads to the
same optimal control as in the original interval case.

\medskip
{\bf Example 3.} Consider the simplest control problem of the form (\ref{1.1})--(\ref{1.3}):
\begin{equation}\label{1.13-5}
y'(t)=u(t), \quad y(0)=2, \quad y(T)=0.
\end{equation}
There are, obviously, a plenty of various controls $u(t)\in L_2(0,T)$ for which such a trajectory~$y(t)$
exists. However, according to Theorems~1 and~2, the optimal control $u(t)$ is unique, while the corresponding
optimal trajectory $y(t)$ can be found by solving the boundary value problem
$$
y''(t)=0, \quad y(0)=2, \quad y(T)=0.
$$
Thus, it can be represented as a straight line in the $(t,y)$-plane connecting the points $(0,2)$ and
$(T,0).$ Assume that we need to find this control $u(t)$ but we only know that $T$ may be equal either to $2$
or to $4.$ The precise information on $T$ will be available starting from $t=1.$

The {\it main question} is: Which control $u(t)$ would be optimal before the time point $t=1?$

\begin{center}
\unitlength=0.7mm
\begin{picture}(132,45)

 \put(8,29){\line(1,0){30}}
 \put(38,29){\line(2,1){27}}
 \put(38,29){\line(4,-1){87.5}}
 \put(8,29){\circle*{1}}
 \put(38,29){\circle*{1}}
 \put(65,42.5){\circle*{1}}
 \put(125.5,7.1){\circle*{1}}
 \put(6,31){\small $0$}
 \put(33,31){\small $1$}
 \put(37,32){\small $0$}
 \put(60.5,43.5){\small $1$}
 \put(38,23){\small $0$}
 \put(122,9.5){\small $3$}
 \put(0,27.5){\small $v_0$}
 \put(31,24){\small $v_1$}
 \put(66,42){\small $v_2$}
 \put(126.4,6){\small $v_3$}
 \put(21,32){\small $e_1$}
 \put(48,39){\small $e_2$}
 \put(82,19.5){\small $e_3$}

\put(42,0){\small Fig. 2. Graph $\Gamma_3$}
\end{picture}
\end{center}

An answer comes when we extend the system (\ref{1.13-5}) to a $3$-star graph $\Gamma_3$ (see Fig.~2) taking
into account both possibilities. In this connection, we denote by $y_1(t)$ a trajectory on the interval
$(0,1),$ which can be only common, and by $\tilde y_2(t)$ and $\tilde y_3(t)$ -- two possible trajectories on
$(1,2)$ and $(1,4),$ respectively. Then we arrive at the corresponding control problem on $\Gamma_3:$
$$
y_j'(t)=u_j(t), \quad 0<t<1, \quad j=1,2, \qquad  y_3'(t)=u_3(t), \quad 0<t<3,
$$
$$
y_1(0)=2, \quad y_1(1)=y_2(0)=y_3(0), \quad y_2(1)=y_3(3)=0,
$$
where $y_j(t)=\tilde y_j(t+1)$ for $j=2,3.$ Assume that the probability of $T=2$ equals $p\in(0,1),$ while
$1-p$ remains to be the probability of $T=4.$ Thus, in the corresponding energy functional (\ref{1.11}), we
will have $\alpha_1=1,$ $\alpha_2=p$ and $\alpha_3=1-p,$ while $T_1=T_2=1$ and $T_3=3.$

According to Theorem~3, the optimal trajectory $[y_1,y_2,y_3]$ solves the boundary value problem
$$
y_j''(t)=0, \quad 0<t<1, \quad j=1,2, \qquad  y_3''(t)=0, \quad 0<t<3,
$$
$$
y_1(0)=2, \quad y_1(1)=y_2(0)=y_3(0), \quad y_1'(1)=py_2'(0)+(1-p)y_3'(0), \quad y_2(1)=y_3(3)=0,
$$
whose solution can be easily checked to have the form
$$
y_1(t)=2(1-t)+\frac{3t}{p+2}, \quad y_2(t)=3\frac{1-t}{p+2}, \quad 0\le t\le1, \qquad y_3(t)=\frac{3-t}{p+2},
\quad 0\le t\le3.
$$
For $p=\frac12,$ this trajectory is given in Fig.~3, where $\tilde y_j(t)=y_j(t-1),\;j=2,3.$ For $j=2,3,$ the
line $y_{1,j}$ shows the optimal trajectory in absence of the scenario corresponding to the edge $e_{5-j}.$

\begin{center}
\unitlength=0.7mm
\begin{picture}(160,87)

 \put(5,15){\vector(1,0){160}}
 \put(10,10){\vector(0,1){73}}
% \put(10,15){\circle*{0.8}}
% \put(10,45){\circle*{0.8}}
 \put(10,75){\circle*{0.8}}
 \put(4,73.3){\small $2$}
 \put(40,15){\circle*{0.8}}
 \put(37.7,9){\small $1(=T_1)$}
 \put(10,60){\circle*{0.8}}
 \put(10,51){\circle*{0.8}}
 \put(10,45){\circle*{0.8}}
 \put(-1,58.5){\small$3/2$}
 \put(-17.3,49.5){\small$\frac3{p+2}=6/5$}
 \put(4.5,43.5){\small$1$}
 \put(70,15){\circle*{0.8}}
 \put(67.7,9){\small $2(=T_1+T_2)$}
 \put(130,15){\circle*{0.8}}
% \put(40,45){\circle*{0.8}}
 \put(40,60){\circle*{0.8}}
 \put(40,51){\circle*{0.8}}
 \put(40,45){\circle*{0.8}}
 \put(10,75){\line(1,-1){60}}
 \put(10,75){\line(2,-1){120}}
 \multiput(40,15)(0,1){61}{\circle*{0.5}}
 \multiput(10,60)(1,0){30}{\circle*{0.5}}
 \multiput(10,51)(1,0){30}{\circle*{0.5}}
 \multiput(10,45)(1,0){30}{\circle*{0.5}}
 \put(10,75){\line(5,-4){30}}
 \put(70,15){\line(-5,6){30}}
 \put(130,15){\line(-5,2){90}}
 \put(127.5,9){\small $4(=T_1+T_3)$}
 \put(33.6,56.8){\small $y_1$}
 \put(49,41){\small $\tilde y_2$}
 \put(55.5,46){\small $\tilde y_3$}
 \put(26.5,47){\small $y_{1,2}$}
 \put(70,47){\small $y_{1,3}$}
 \put(164,17){\small $t$}
 \put(10,84){\small $y$}

\put(11,-3){\small Fig. 3. The optimal trajectories in Example~3 for $p=\frac12$}
\end{picture}
\end{center}

\medskip
Although Fig.~3 shows that the trajectory $y_{1,2}$ is shorter than the composite trajectory $[y_1,\tilde
y_2],$ while $y_{1,3}$ is shorter than $[y_1,\tilde y_3],$ both composite trajectories allow one to take more
advantageous position at the time point $t=1,$ before which the further scenario is uncertain.

This is like when someone is going on a trip from $v_0$ to $v_1$ and then to either $v_2$ or $v_3.$ But the
choice between $v_2$ and $v_3$ can be made only in $v_1.$ For this reason, the travel bag should be assembled
in such a way as to minimize, for example, the weight of potentially useless items that could not be left in
$v_1$ and the cost of potentially missing items that should be purchased in it.

Finally, recall that $p$ equals the probability of the scenario corresponding to the edge~$e_2.$ Obviously,
the point $y_1(1)=\frac3{p+2},$ which is indicated in Fig.~4 for $p=\frac12,$ sweeps the interval
$(y_{1,2}(1),y_{1,3}(1))=(1,\frac32)$ as soon as $p$ ranges over $(0,1).$ In particular, the larger $p$ is,
the closer the composite trajectory $[y_1,\tilde y_2]$ is to the line $y_{1,2},$ while $[y_1,\tilde y_3]$
tends to $y_{1,3}$ when $p\to0.$

\medskip
In the next section, we generalize the above variational problem to an arbitrary tree with a finite number of
edges. Although the tree can possess originally infinite boundary edges as in $\Gamma_m,$ we will actually
deal with a compact tree ${\cal T}$ obtained by cutting the infinite edges $e_j$ at the points~$T_j$ starting
from which the control system should be damped. For definiteness, we will keep assuming that the delay
parameter $\tau$ is less than all $T_j,$ i.e. less than the length of each edge in ${\cal T}.$ However,
applying the settings from Section~7 of \cite{But23}, one can consider a more general case, e.g., when
$2\tau<T,$ where $T$ is now the height of a tree.
\\

{\bf 3. Statement of the variational problem on a tree}
\\

Consider a compact rooted tree ${\cal T}$ with the set of vertices $\{v_0,v_1,\ldots, v_m\}$ and the set of edges $\{e_1,\ldots, e_m\}.$ Let
$\{v_0,v_{d+1},\ldots, v_m\}$ be boundary vertices, i.e. each of them is incident to only one (boundary) edge. The remaining vertices
$\{v_1,\ldots, v_d\}$ are internal.

Without loss of generality, we agree that each edge $e_j,\;j=\overline{1,m},$ emanates from the corresponding
vertex $v_{k_j}$ and terminates at $v_j$ and write $e_j=[v_{k_j},v_j],$ where $k_1=0.$ The vertex $v_0$ will
be labelled as {\it root}. For example, in Fig.~4, we have $m=9$ and $d=3,$ while
$$
k_1=0, \quad k_2=k_3=1, \quad k_4=k_5=2, \quad k_6=k_7=k_8=k_9=3.
$$

\begin{center}
\unitlength=0.7mm
\begin{picture}(180,100)

 \put(90,18){\circle*{1}}         \put(87,13){\small $v_0$}
 \put(90,46){\circle*{1}}         \put(87,50){\small $v_1$}
 \put(66,61.95){\circle*{1}}      \put(58,64){\small $v_2$}
 \put(114,61.95){\circle*{1}}     \put(116,64){\small $v_3$}
 \put(40,49){\circle*{1}}         \put(33,46){\small $v_4$}
 \put(66,89.95){\circle*{1}}      \put(63,92.5){\small $v_5$}
 \put(94,82){\circle*{1}}         \put(87,84){\small $v_6$}
 \put(123,89){\circle*{1}}        \put(121,91.5){\small $v_7$}
 \put(142,61.95){\circle*{1}}     \put(143,61){\small $v_8$}
 \put(123,35){\circle*{1}}        \put(121.5,30.5){\small $v_9$}

 \put(90,18){\line(0,1){28}}      \put(83,32){\small $e_1$}
 \put(90,46){\line(-3,2){24}}     \put(72,50){\small $e_2$}
 \put(90,46){\line(3,2){24}}      \put(100,50){\small $e_3$}
 \put(66,61.95){\line(-2,-1){26}} \put(51.5,51.5){\small $e_4$}
 \put(66,61.95){\line(0,1){28}}   \put(66,76){\small $e_5$}
 \put(114,61.95){\line(-1,1){20}} \put(97,70){\small $e_6$}
 \put(114,61.95){\line(1,3){9}}   \put(112,78){\small $e_7$}
 \put(114,61.95){\line(1,0){28}}  \put(128,64){\small $e_8$}
 \put(114,61.95){\line(1,-3){9}}  \put(119,47){\small $e_9$}

 \put (58,0){\small Fig. 4. A non-star tree}

\end{picture}
\end{center}

Thus, $k_j$ generates some (non-injective when $m>d+1)$ mapping of the set $\{1,\ldots,m\}$ onto $\{0,1,\ldots,d\},$ which uniquely
determines the structure of ${\cal T}.$ Specifically, for each $j=\overline{0,d},$ the set $\{e_\nu\}_{\nu\in V_j},$ where
\begin{equation}\label{2.0}
V_j:=\{\nu:k_\nu=j\},
\end{equation}
coincides with the set of edges emanating from the vertex $v_j.$ In particular,
% $\#\{e_\nu\}_{\nu\in V_0}=1$
$\#V_0=1$ since $v_0$ is a boundary vertex. Put $k_j^{<0>}:=j$ and $k_j^{<\nu+1>}:=k_{k_j^{<\nu>}}$ for
$\nu=\overline{0,\nu_j},$ where $\nu_j$ is determined by $k_j^{<\nu_j>}=1.$ Then, for each
$j=\overline{1,m},$ the chain of edges $\{e_{k_j^{<\nu>}}\}_{\nu=\overline{0,\nu_j}}$ forms the unique simple
path between the vertex $v_j$ and the root. Denote the length of the edge $e_j$ by $T_j.$ Then the value
$T:=\max_{j=\overline{d+1,m}}\sum_{\nu=0}^{\nu_j}T_{k_j^{<\nu>}}$ is called {\it height} of the tree ${\cal
T}.$

Let each edge $e_j$ be parametrized by the variable $t\in[0,T_j]$ so that $t=0$ and $t=T_j$ correspond to its
ends $v_{k_j}$ and $v_j,$ respectively. By a function $y$ on ${\cal T},$ we mean an $m\!$-tuple
$y=[y_1,\ldots,y_m]$ whose component $y_j$ is defined on the edge $e_j,$ i.e. $y_j=y_j(t),\,t\in[0,T_j].$ We
also fix $\tau\ge0$ and say that the function $y$ is defined on the extended tree ${\cal T}_\tau$ if it is
defined on ${\cal T}$ and its first component $y_1(t)$ is defined also for $t\in[-\tau,0).$

For definiteness, let $\tau<T_j$ for all $j=\overline{1,m}.$ Consider the following Cauchy problem on~${\cal T}_\tau:$
\begin{equation}\label{2.1}
\ell_j y(t):=y_j'(t)+b_jy_j(t)+c_jy_j(t-\tau)=u_j(t), \quad 0<t<T_j, \quad j=\overline{1,m},
\end{equation}
\begin{equation}\label{2.2}
y_j(t)=y_{k_j}(t+T_{k_j}), \quad t\in(-\tau,0), \quad j=\overline{2,m},
\end{equation}
\begin{equation}\label{2.2-1}
y_j(0)=y_{k_j}(T_{k_j}), \quad j=\overline{2,m},
\end{equation}
\begin{equation}\label{2.3}
y_1(t)=\varphi(t)\in W_2^1[-\tau,0], \quad t\in[-\tau,0],
\end{equation}
where $b_j,c_j\in{\mathbb R}$ and $u_j\in L_2(0,T_j)$ for $j=\overline{1,m}.$ While the $j$th equation in (\ref{2.1}) is defined on the edge
$e_j$ of ${\cal T},$ relations (\ref{2.2-1}) become matching conditions at the internal vertices. Relations (\ref{2.2}) are initial-function
conditions for all equations in (\ref{2.1}) except the first one. They mean that the delay propagates through all internal vertices.
Condition (\ref{2.3}) determines the prehistory of the process for the entire tree, where the function $\varphi(t)$ is real valued and known.

We note that the problem (\ref{1.5})--(\ref{1.8}), obviously, coincides with the problem
(\ref{2.1})--(\ref{2.3}) for $d=1.$ As in the preceding section, one can show that the latter problem also
has a unique solution
$$
[y_1,\ldots,y_m]\in \bigoplus_{j=1}^m W_2^1[0,T_j].
$$

\medskip
{\bf Remark 1.} Although conditions (\ref{2.2}) and (\ref{2.2-1}) can be formally combined by allowing $t$ in (\ref{2.2}) to take the zero
value, we prefer not to do it since these two conditions are of a different nature. Moreover, one can consider their more general forms
$$
a_jy_j(t)=y_{k_j}(t+T_{k_j}), \;\; t\in(-\tau,0), \quad \tilde a_jy_j(0)=y_{k_j}(T_{k_j}), \quad j=\overline{2,m},
$$
where $a_j$ and $\tilde a_j$ may differ. Meanwhile, it is natural to combine (\ref{2.2}) with (\ref{2.1}):
$$
\ell_j y(t)=y_j'(t)+b_jy_j(t)+c_j\left\{\begin{array}{cc}
y_{k_j}(t-\tau+T_{k_j}), & 0<t<\tau,\\[3mm]
y_j(t-\tau), & \tau<t<T_j,
\end{array}\right. \quad j=\overline{2,m}.
$$
However, the separation of (\ref{2.1}) and (\ref{2.2}) is more convenient and shows the succession to (\ref{1.1}). Analogously,
condition (\ref{2.3}) also can be naturally split into the two ideologically different ones:
$$
y_1(t)=\varphi(t), \;\; t\in[-\tau,0), \quad {\rm and} \quad y_1(0)=y_0, \quad y_0:=\varphi(0),
$$
where the first condition initializes the functional part $c_1y_1(t-\tau)$ of the expression $\ell_1y(t),$
while the second one becomes an initial condition for its differential part $y_1'(t).$ In general, one can
consider the situation when $\varphi(0)\ne y_0.$ However, for simplicity and for consistency with
\cite{Skub-94}, we keep $\varphi(0)=y_0$ and use here the combined form (\ref{2.3}) as more concise.

\medskip
Analogously to the preceding section, we intend to find a control function
$$
u=[u_1,\ldots,u_m]\in L_2({\cal T}):=\bigoplus_{j=1}^m L_2(0,T_j)
$$
that leads to the equilibrium state
\begin{equation}\label{2.4}
y_j(t)=0, \quad t\in[T_j-\tau,T_j], \quad j=\overline{d+1,m},
\end{equation}
and minimizes all $\|u_j\|_{L_2(0,T_j)}^2$ with some positive weights $\alpha_j,\;j=\overline{1,m}.$

Thus, we arrive at the variational problem:
\begin{equation}\label{2.5}
{\cal J}(y):=\sum_{j=1}^m \alpha_j\int_0^{T_j}(\ell_j y(t) )^2\,dt\to\min
\end{equation}
for the functions $y=[y_1,\ldots,y_m]$ defined on ${\cal T}_\tau$ under the conditions (\ref{2.2})--(\ref{2.4}).

For brevity, we introduce the designation $\ell y:=[\ell_1y,\ldots,\ell_my]$ and agree that taking ${\cal
J}(y)$ and $\ell y$ as well as $\ell_j y$ for $j=\overline{2,m}$ of any function $y$ on ${\cal T}$
automatically means the application of the initial-function conditions (\ref{2.2}).
\\

{\bf 4. Reduction to a boundary value problem}
\\

Consider the real Hilbert space $W_2^k({\cal T}_\tau):=W_2^k[-\tau,T_1]\oplus\bigoplus_{j=2}^m W_2^k[0,T_j]$ with the natural inner product
$$
(y,z)_{W_2^k({\cal T}_\tau)}:=(y_1,z_1)_{W_2^k[-\tau,T_1]}+\sum_{j=2}^m (y_j,z_j)_{W_2^k[0,T_j]},
$$
where $y=[y_1,\ldots,y_m]$ and $z=[z_1,\ldots,z_m],$ while $(f,g)_{W_2^k[a,b]}=\sum_{\nu=0}^k(f^{(\nu)},g^{(\nu)})_{L_2(a,b)}$ is the inner
product in $W_2^k[a,b]$ and $(\,\cdot\,,\,\cdot\,)_{L_2(a,b)}$ is the one in $L_2(a,b).$ In particular, $W_2^0({\cal T}_0)=L_2({\cal T}).$

Denote by ${\cal W}$ the closed subspace of $W_2^1({\cal T}_\tau)$ consisting of all $m\!$-tuples $[y_1,\ldots,y_m]\in W_2^1({\cal T}_\tau)$
that obey the matching conditions (\ref{2.2-1}), the target conditions (\ref{2.4}), and  $y_1(t)=0$ on $[-\tau,0].$

Obviously, ${\cal W}$ is a subspace also of $W_2^1({\cal T}),$ as well as of $W_2^1(\widetilde{\cal T}),$ where $W_2^k({\cal T}):=W_2^k({\cal
T}_0)$ while $W_2^k(\widetilde{\cal T}):=\bigoplus_{j=1}^d W_2^k[0,T_j]\oplus\bigoplus_{j=d+1}^m W_2^k[0,T_j-\tau].$ In other words,
$W_2^k(\widetilde{\cal T})$ differs from $W_2^k({\cal T})$ only by replacing $T_j$ with $T_j-\tau$ for $j=\overline{d+1,m}.$

As in the case of an interval, the optimal trajectory on a tree possesses additional smoothness (see Theorem~5 at the end of this section).
Specifically, besides $W_2^1({\cal T}_\tau),$ it belongs also to the set $W_2^2(\widetilde{\cal T})$ introduced above. Moreover, as in the
case of a star-shaped graph, the optimal trajectory obeys Kirchhoff-type conditions at the internal vertices (\ref{3.14}), which, in turn,
generalize the continuity of the first derivative at each internal point of an edge.

\medskip
{\bf Lemma 1. }{\it If $y\in W_2^1({\cal T}_\tau)$ is a solution of the variational problem (\ref{2.2})--(\ref{2.5}), then
\begin{equation}\label{3.1}
B(y,w):=\sum_{j=1}^m\alpha_j\int_0^{T_j} \ell_j y(t)\ell_j w(t)\,dt=0 \quad \forall \; w\in {\cal W}.
\end{equation}
Conversely, if $y\in W_2^1({\cal T}_\tau)$ obeys (\ref{2.2-1})--(\ref{2.4}) and (\ref{3.1}), then $y$
is a solution of (\ref{2.2})--(\ref{2.5}).}

\medskip
{\it Proof.} Let $y\in W_2^1({\cal T}_\tau)$ be a solution of (\ref{2.2})--(\ref{2.5}). Then for any $w\in {\cal W},$ the sum $y+sw$ belongs
to $W_2^1({\cal T}_\tau)$ whenever $s\in{\mathbb R},$ and obeys conditions (\ref{2.2-1})--(\ref{2.4}). Put
$$
F(s):={\cal J}(y+sw)={\cal J}(y)+2sB(y,w)+s^2{\cal J}(w).
$$
Since ${\cal J}(y+sw)\ge {\cal J}(y)$ for all $s\in{\mathbb R},$ we have $0=F'(0)=2B(y,w).$

Conversely, for any $y\in W_2^1({\cal T}_\tau)$ obeying (\ref{2.2-1})--(\ref{2.4}), the fulfilment of (\ref{3.1})
implies
$$
{\cal J}(y+w)={\cal J}(y)+2B(y,w)+{\cal J}(w)\ge {\cal J}(y)
$$
for all $w\in {\cal W},$ which gives (\ref{2.5}).  $\hfill\Box$

\medskip
Further, applying (\ref{2.2}) to $w=[w_1,\ldots,w_m]\in{\cal W},$ one can represent
$$
\int_0^{T_j}\ell_jy(t)w_j(t-\tau)\,dt= \int_0^{T_j-\tau}\ell_jy(t+\tau)w_j(t)\,dt \qquad\qquad\qquad\qquad\quad
$$
$$
\;\qquad\qquad\qquad\qquad +\int_{T_{k_j}-\tau}^{T_{k_j}}\ell_jy(t+\tau-T_{k_j})w_{k_j}(t)\,dt, \quad j=\overline{1,m},
$$
where $T_0=0$ and $w_0=0.$ According to (\ref{2.0}), we have the summation rule
\begin{equation}\label{3.1-1}
\sum_{j=2}^m A_j=\sum_{j=1}^d\sum_{\nu\in V_j}A_\nu \quad ({\rm for \;\; any \;\; values} \;\; A_2,\ldots,A_m).
\end{equation}
Hence, multiplying the preceding relation with $\alpha_jc_j$ and then summing up, we obtain
$$
\sum_{j=1}^m \alpha_jc_j\int_0^{T_j}\ell_jy(t)w_j(t-\tau)\,dt= \sum_{j=1}^m
\alpha_jc_j\int_0^{T_j-\tau}\ell_jy(t+\tau)w_j(t)\,dt\qquad\qquad\qquad\qquad
$$
$$
\qquad\qquad\qquad\qquad +\sum_{j=1}^d \sum_{\nu\in V_j} \alpha_\nu c_\nu \int_{T_j-\tau}^{T_j}\ell_\nu y(t+\tau-T_j)w_j(t)\,dt.
$$
Thus, in accordance with the definition in (\ref{2.1}), one can rewrite (\ref{3.1}) in the equivalent form
\begin{equation}\label{3.2}
B(y,w)=\sum_{j=1}^m \Big(\alpha_j\int_0^{T_j} \ell_j y(t) w_j'(t)\,dt +\int_0^{T_j} (\alpha_jb_j\ell_j y(t) +\tilde\ell_j y(t))w_j(t)
\,dt\Big)=0, \quad w\in {\cal W},
\end{equation}
where
\begin{equation}\label{3.3}
\tilde\ell_j y(t):=\left\{\begin{array}{cc}
\alpha_jc_j\ell_j y(t+\tau), & 0<t<T_j-\tau, \;\; j=\overline{1,m},\\[3mm]
\displaystyle\sum_{\nu\in V_j}\alpha_\nu c_\nu\ell_\nu y(t+\tau-T_j), & T_j-\tau< t<T_j, \;\; j=\overline{1,d},
\end{array}\right.
\end{equation}
while for $T_j-\tau<t< T_j$ and $j=\overline {d+1,m},$ the expression $\tilde\ell_jy(t)$ can be defined as zero.

\medskip
{\bf Lemma 2. }{\it Let for some $g_j(t),\,f_j(t)\in L_2(0,l_j),\;j=\overline{1,m},$ the relation
\begin{equation}\label{3.4}
\sum_{j=1}^m \Big(\int_0^{l_j} g_j(t) w_j'(t)\,dt +\int_0^{l_j} f_j(t)w_j(t) \,dt\Big)=0
\end{equation}
be fulfilled with each $m\!$-tuple of $w_j(t)\in W_2^1[0,l_j],\;j=\overline{1,m},$ obeying the conditions
\begin{equation}\label{3.5}
w_1(0)=0, \quad w_j(0)=w_{k_j}(l_{k_j}), \;\; j=\overline{2,m}, \quad w_j(l_j)=0, \;\;  j=\overline{d+1,m}.
\end{equation}

Then $g_j(t)\in W_2^1[0,l_j]$ for $j=\overline{1,m}$ and}
\begin{equation}\label{3.6}
g_j(l_j)=\sum_{\nu\in V_j}g_\nu(0), \quad j=\overline{1,d}.
\end{equation}

\medskip
{\it Proof.} Choose antiderivatives $F_j(t)\in W_2^1[0,l_j]$ of $f_j(t),$ i.e. $F_j'(t)=f_j(t),$ such that
\begin{equation}\label{3.7}
F_j(l_j)=\sum_{\nu\in V_j} F_\nu(0), \quad j=\overline{1,d},
\end{equation}
which is possible. Indeed, for $j=\overline{d+1,m},$ one can choose any antiderivative $F_j(t)$ of $f_j(t),$ while the remaining $F_j(t),$
$j=\overline{1,d},$ are recurrently defined by conditions (\ref{3.7}).

Then, integrating by parts in (\ref{3.4}), we obtain
\begin{equation}\label{3.8}
\sum_{j=1}^m \int_0^{l_j} (g_j(t)-F_j(t)) w_j'(t)\,dt +\sum_{j=1}^m (F_j(t)w_j(t))\Big|_{t=0}^{l_j}=0.
\end{equation}
In accordance with (\ref{2.0}), the second relation in (\ref{3.5}) gives $w_j(l_j)=w_\nu(0)$ for $\nu\in V_j$
and $j=\overline{1,d}.$ Thus, using (\ref{3.1-1}), (\ref{3.5}) and (\ref{3.7}), we obtain
\begin{equation}\label{3.9}
\begin{array}{c}
\displaystyle\sum_{j=1}^m (F_j(t)w_j(t))\Big|_{t=0}^{l_j}= \sum_{j=1}^d \Big(F_j(l_j)w_j(l_j)
-\sum_{\nu\in V_j} F_\nu(0)w_\nu(0)\Big) \qquad\qquad\qquad\\[3mm]
\displaystyle \qquad\qquad\qquad\qquad\qquad\qquad= \sum_{j=1}^d w_j(l_j)\Big(F_j(l_j)-\sum_{\nu\in V_j} F_\nu(0)\Big)=0.
\end{array}
\end{equation}
Hence, for each $j=\overline{1,m},$ after putting $w_k(t)\equiv0$ for all $k\ne j,$ relation (\ref{3.8}) takes the form
$$
\int_0^{l_j} (g_j(t)-F_j(t)) w_j'(t)\,dt =0, \quad j=\overline{1,m},
$$
for all functions $w_j(t)\in W_2^1[0,l_j]$ such that $w_j(0)=w_j(l_j)=0.$ Since the orthogonal complement in $L_2(0,l_j)$ of their
derivatives $w_j'(t)$ consists of constants, we obtain
\begin{equation}\label{3.10}
g_j(t)-F_j(t)=C_j \equiv const, \quad j=\overline{1,m}.
\end{equation}
This gives the first assertion of the lemma, namely that $g_j(t)\in W_2^1[0,l_j]$ for $j=\overline{1,m}.$

Further, substituting (\ref{3.10}) into (\ref{3.8}) and using (\ref{3.9}) along with (\ref{3.1-1}) and
(\ref{3.5}), we obtain
$$
0= \sum_{j=1}^m C_j \int_0^{l_j} w_j'(t)\,dt =\sum_{j=1}^d w_j(l_j)\Big(C_j-\sum_{\nu\in V_j} C_\nu\Big).
$$
Hence, by virtue of the arbitrariness of $w_j(l_j)$ for $j=\overline{1,d},$ we have
\begin{equation}\label{3.11}
C_j=\sum_{\nu\in V_j} C_\nu, \quad j=\overline{1,d}.
\end{equation}
Finally, substituting (\ref{3.10}) into (\ref{3.7}), we arrive at
$$
g_j(l_j)-\sum_{\nu\in V_j} g_\nu(0)=C_j-\sum_{\nu\in V_j} C_\nu, \quad j=\overline{1,d},
$$
which along with (\ref{3.11}) gives (\ref{3.6}). $\hfill\Box$

\medskip
Denote by ${\cal B}$ the boundary value problem for the second-order functional-differential equations
\begin{equation}\label{3.17}
{\cal L}_jy(t):=-\alpha_j(\ell_j y)'(t)+\alpha_jb_j\ell_j y(t) +\tilde\ell_j y(t)=0, \quad 0<t< l_j, \quad j=\overline{1,m},
\end{equation}
under the conditions (\ref{2.2})--(\ref{2.4}) along with the conditions
\begin{equation}\label{3.14}
\alpha_jy_j'(T_j) +\beta_j y_j(T_j) +\gamma_j y_j(T_j-\tau)=\sum_{\nu\in V_j}\alpha_\nu y_\nu'(0), \quad j=\overline{1,d},
\end{equation}
where the expressions $\tilde\ell_j y(t)$ are defined in (\ref{3.3}) and
\begin{equation}\label{3.12}
l_j:=T_j,\quad j=\overline{1,d}, \qquad l_j:=T_j-\tau, \quad j=\overline{d+1,m},
\end{equation}
\begin{equation}\label{3.15}
\beta_j :=\alpha_jb_j -\sum_{\nu\in V_j}\alpha_\nu b_\nu, \quad \gamma_j :=\alpha_jc_j -\sum_{\nu\in V_j}\alpha_\nu c_\nu, \quad
j=\overline{1,d}.
\end{equation}

The following lemma holds.

\medskip
{\bf Lemma 3. }{\it If $y\in W_2^1({\cal T}_\tau)$ obeys conditions (\ref{2.2-1})--(\ref{2.4}) and
(\ref{3.1}), then $y\in W_2^2(\widetilde{\cal T})$ and it solves the boundary value problem ${\cal B}.$
Conversely, any solution $y$ of ${\cal B}$ obeys (\ref{3.1}).}

\medskip
{\it Proof.} Taking into account that (\ref{3.1}) is equivalent to (\ref{3.2}) and applying Lemma~2 to (\ref{3.2})
under the settings (\ref{3.12}), we get $\ell_j y(t)\in W_2^1[0,l_j]$ for $j=\overline{1,m}$ and
\begin{equation}\label{3.12-0}
\alpha_j\ell_jy(T_j)=\sum_{\nu\in V_j}\alpha_\nu\ell_\nu y(0), \quad j=\overline{1,d}.
\end{equation}
Recalling the definition of $\ell_jy$ in (\ref{2.1}) and using (\ref{2.2})--(\ref{2.3}) as well as (\ref{3.12}), we obtain the inclusion
$y_j'(t)=\ell_j y(t)-b_jy_j(t)-c_jy_j(t-\tau) \in W_2^1[0,l_j]$ for $j=\overline{1,m},$ i.e. $y\in W_2^2(\widetilde{\cal T}).$ Hence, we can
rewrite (\ref{3.12-0}) in the equivalent form
\begin{equation}\label{3.12-1}
y_j'(T_j)+b_jy_j(T_j)+c_jy_j(T_j-\tau) =\frac1{\alpha_j}\sum_{\nu\in V_j}\alpha_\nu(y_\nu'(0)+b_\nu y_\nu(0)+c_\nu y_\nu(-\tau)), \quad
j=\overline{1,d},
\end{equation}
where, by virtue of (\ref{2.2}), the right-hand limits
$$
y_\nu(-\tau):=\lim_{t\to(-\tau)^+}y_\nu(t) =\lim_{t\to(-\tau)^+}y_{k_\nu}(t+T_{k_\nu})
=y_{k_\nu}(T_{k_\nu}-\tau),
$$
obviously, exist. Thus, relations (\ref{3.12-1}) along with (\ref{2.0}), (\ref{2.2}) and (\ref{2.2-1}) give (\ref{3.14}) with (\ref{3.15}).

Finally, integrating by parts in (\ref{3.2}) and acting as in (\ref{3.9}), we arrive at
\begin{equation}\label{3.16}
B(y,w)=\sum_{j=1}^d w_j(l_j)\Big(\alpha_j\ell_jy(l_j)-\sum_{\nu\in V_j} \alpha_\nu\ell_\nu y(0)\Big) +\sum_{j=1}^m \int_0^{l_j}{\cal L}_jy(t)
w_j(t) \,dt=0,
\end{equation}
which, by virtue of (\ref{3.12-0}) and the variety of $w_j(t),$ gives (\ref{3.17}).

Conversely, let $y$ be a solution of the problem ${\cal B}.$ Since (\ref{3.14}) is equivalent to (\ref{3.12-0}), the second equality in
(\ref{3.16}) holds, which gives (\ref{3.1}).  $\hfill\Box$

\medskip
{\bf Remark 2.} We also established that any solution of the variational problem (\ref{2.2})--(\ref{2.5})
obeys nonlocal Kirchhoff-type conditions (\ref{3.14}) at the internal vertices. In particular, if
$\alpha_\nu=\alpha_j,$ $\nu\in V_j,$ and $\beta_j=\gamma_j=0$ for certain $j\in\{1,\ldots,d\},$ then the
corresponding relation in (\ref{3.14}) becomes a usual Kirchhoff condition, which complements the standing
continuity conditions $y_\nu(0)=y_j(T_j)$ for $\nu\in V_j$ to the so-called standard matching conditions at
the internal vertex~$v_j.$

\medskip

{\bf Remark 3.} By virtue of the left-hand equalities in (\ref{3.1}) and (\ref{3.16}) along with the equivalence of (\ref{3.14}) and
(\ref{3.12-0}) under the continuity conditions (\ref{2.2-1}), the problem ${\cal B}$ is self-adjoint. Specifically, this means that the
relation $({\cal L}y,z)_{L_2(\widetilde{\cal T})}=(y,{\cal L}z)_{L_2(\widetilde{\cal T})}$ is fulfilled for all $y,z\in{\cal D},$ where
$$
{\cal D}:=\{y: y\in{\cal W}\cap W_2^2(\widetilde{\cal T})\; {\rm and}\; y \;{\rm obeys\; conditions}\; (\ref{3.14})\}
$$
and ${\cal L}y:=[{\cal L}_1y,\dots,{\cal L}_my],$ while $(\,\cdot\,,\,\cdot\,)_{L_2(\widetilde{\cal T})}$ is the inner product in
$L_2(\widetilde{\cal T}):=W_2^0(\widetilde{\cal T}).$ Moreover, one can show that the operator ${\cal L}$ with domain ${\cal D}$ is a
self-adjoint operator in $L_2(\widetilde{\cal T}),$ i.e. the set ${\cal D}$ is dense in $L_2(\widetilde{\cal T})$ and the domain of the
adjoint operator ${\cal L}^*$ is not wider than ${\cal D}.$ Lemma~6 in the next section also implies that this operator is positive definite.

\medskip
Combining Lemma~1 and Lemma~3, we arrive at the main result of this section.

\medskip
{\bf Theorem 5. }{\it A function $y\in W_2^1({\cal T}_\tau)$ is a solution of the variational problem (\ref{2.2})--(\ref{2.5}) if and only if
it belongs to $W_2^2(\widetilde{\cal T})$ and solves the boundary value problem ${\cal B}.$}
\\

{\bf 5. The unique solvability}
\\

In this section, we establish the unique solvability of the boundary value problem ${\cal B}$ and thus, according to Theorem~5, of the
variational problem (\ref{2.2})--(\ref{2.5}).

We begin with the following three auxiliary assertions.

\medskip
{\bf Lemma 4. }{\it There exists $C$ such that}
\begin{equation}\label{4.0}
\|\ell w\|_{L_2({\cal T})}^2\le C\|w\|_{W_2^1(\widetilde{\cal T})}^2 \quad \forall w\in {\cal W}.
\end{equation}

\medskip
{\it Proof.} Using the definition in (\ref{2.1}) and the inequality
\begin{equation}\label{4.3-1}
(a_1+\ldots+a_n)^2\le n(a_1^2+\ldots+a_n^2), \quad a_1,\ldots,a_n\in{\mathbb R},
\end{equation}
for $n=3,$ we obtain
\begin{equation}\label{4.00}
\|\ell w\|_{L_2({\cal T})}^2\le 3\sum_{j=1}^m \int_0^{T_j}(w_j'(t))^2\,dt +3\sum_{j=1}^m b_j^2\int_0^{T_j}w_j^2(t)\,dt +3\sum_{j=1}^m
c_j^2\int_0^{T_j}w_j^2(t-\tau)\,dt.
\end{equation}
Then, applying (\ref{2.2}) to $w,$ we calculate
\begin{equation}\label{4.01}
\int_0^{T_j}w_j^2(t-\tau)\,dt=\|w_j\|_{L_2(0,T_j-\tau)}^2 + \|w_{k_j}\|_{L_2(T_{k_j}-\tau,T_{k_j})}^2, \quad j=\overline{1,m},
\end{equation}
where $T_0=0$ and $w_0=0,$ which along with (\ref{4.00}) gives (\ref{4.0}) with $C$ independent of $w.$ $\hfill\Box$

\medskip
{\bf Lemma 5. }{\it There exists $c>0$ such that
\begin{equation}\label{4.1}
\|w'\|_{L_2({\cal T})}^2\ge c\|w\|_{W_2^1(\widetilde{\cal T})}^2 \quad \forall w\in {\cal W},
\end{equation}
where $w'=[w_1',\ldots, w_m'].$}

\medskip
{\it Proof.} Fix $j\in\{1,\ldots,m\}$ and let $\{e_{j_\nu}\}_{\nu=\overline{0,s}}$ be the simple path between $v_j$ and $v_0,$ i.e.
$j_\nu=k_j^{<\nu>}$ for $\nu=\overline{0,s},$ and $k_j^{<s>}=1.$ Since any function $w\in {\cal W}$ obeys the continuity conditions
(\ref{2.2-1}) as well as the boundary condition $w_1(0)=0,$ we have the representations
$$
w_{j_\nu}(t)=w_{j_{\nu+1}}(T_{j_{\nu+1}}) +\int_0^t w_{j_\nu}'(\xi)\,d\xi, \;\; \nu=\overline{0,s-1}, \quad w_1(t)= \int_0^t w_1'(\xi)\,d\xi.
$$
Successively using them and taking into account that $j_0=j$ and $j_s=1,$ we obtain
$$
w_j(t)=\sum_{\nu=1}^s\int_0^{T_{j_\nu}} w_{j_\nu}'(\xi)\,d\xi +\int_0^t w_j'(\xi)\,d\xi.
$$
By virtue of inequality (\ref{4.3-1}), there exists $C_j,$ which depends only on $T_{j_\nu}$ for $\nu=\overline{0,s},$ such that
$$
\|w_j\|_{L_2(0,T_j)}^2\le C_j\sum_{\nu=0}^s \|w_{j_\nu}'\|_{L_2(0,T_{j_\nu})}^2.
$$
Thus, we arrive at the estimate
$$
\|w\|_{L_2({\cal T})}^2 \le C_0\|w'\|_{L_2({\cal T})}^2,
$$
where $C_0$ depends only on $k_j$ and $T_j$ for $j=\overline{1,m}.$ Since
$$
\|w\|_{W_2^1({\cal T})}^2=\|w\|_{L_2({\cal T})}^2 +\|w'\|_{L_2({\cal T})}^2,
$$
the latter estimate gives (\ref{4.1}) with $c=(1+C_0)^{-1}.$ $\hfill\Box$

\medskip
{\bf Lemma 6. }{\it There exists $C_1>0$ such that}
\begin{equation}\label{4.2}
{\cal J}(w)\ge C_1\|w\|_{W_2^1(\widetilde{\cal T})}^2 \quad \forall w\in {\cal W}.
\end{equation}

\medskip
{\it Proof.} To the contrary, let there exist $w_{(n)}\in {\cal W},$ $n\in{\mathbb N},$ such that
$$
{\cal J}(w_{(n)})\le \frac1n\|w_{(n)}\|_{W_2^1(\widetilde{\cal T})}^2.
$$
Assuming without loss of generality that $\|w_{(n)}\|_{W_2^1(\widetilde{\cal T})}=1,$ we arrive at the inequalities
\begin{equation}\label{4.3}
{\cal J}(w_{(n)})\le \frac1n, \quad n\in{\mathbb N}.
\end{equation}
Using the definition of $\ell_jy$ and inequality (\ref{4.3-1}) for $n=3,$ we obtain
$$
(w_j'(t))^2 \le3((\ell_j w(t))^2 +b_j^2w_j^2(t) +c_j^2w_j^2(t-\tau)), \quad j=\overline{1,m}.
$$
Integrate this inequality from $0$ to $T_j$ and multiply with $\alpha_j.$ Then summing up with respect to $j$
and recalling the definition in (\ref{2.5}), we arrive at
$$
 \frac13\sum_{j=1}^m \alpha_j \int_0^{T_j}
(w_j'(t))^2\,dt\le {\cal J}(w)+ \sum_{j=1}^m \alpha_j \int_0^{T_j}(b_j^2w_j^2(t) + c_j^2w_j^2(t-\tau))\,dt,
$$
which along with (\ref{4.01}) and Lemma~5 implies
\begin{equation}\label{4.5}
\frac{\alpha c}3\|w\|_{W_2^1(\widetilde{\cal T})}^2 \le {\cal J}(w)  +K\|w\|_{L_2({\cal T})}^2, \quad
\alpha:=\min_{j=\overline{1,m}}\alpha_j>0.
\end{equation}

Further, by virtue of the compactness of the embedding operator from $W_2^1({\cal T})$ into $L_2({\cal T}),$ there exists a sequence
$\{w_{(n_k)}\}_{k\in{\mathbb N}}$ that converges in $L_2({\cal T}).$ Inequality (\ref{4.5}) gives
$$
\frac{\alpha c}3\|w_{(n_k)}-w_{(n_l)}\|_{W_2^1(\widetilde{\cal T})}^2 \le{\cal J}(w_{(n_k)}-w_{(n_l)}) +
K\|w_{(n_k)}-w_{(n_l)}\|_{L_2({\cal T})}^2.
$$
Moreover, using (\ref{4.3-1}) for $n=2$ along with (\ref{4.3}), we get
$$
{\cal J}(w_{(n_k)}-w_{(n_l)})\le\frac2{n_k}+\frac2{n_l}.
$$
Thus, $\{w_{(n_k)}\}_{k\in{\mathbb N}}$ is a Cauchy sequence also in ${\cal W}.$ Let $w_{(0)}$ be its limit therein.

By virtue of Lemma~4, the convergence of $w_{(n_k)}$ to $w_{(0)}$ in ${\cal W}$ implies  $\ell w_{(n_k)} \to \ell w_{(0)}$ in $L_2({\cal
T}).$ Then, due to (\ref{4.3}), we have
$$
\|\ell w_{(0)}\|_{L_2({\cal T})}^2= \lim_{k\to\infty}\|\ell w_{(n_k)}\|_{L_2({\cal T})}^2\le\frac1\alpha\lim_{k\to\infty}{\cal
J}(w_{(n_k)})=0,
$$
i.e. $\ell w_{(0)}=0.$ Thus, $w_{(0)}$ solves the Cauchy problem (\ref{2.1})--(\ref{2.3}) with $\varphi(t)\equiv 0$ and $u_j(t)\equiv 0$ for
$j=\overline{1,m}.$ Hence, we have $w_{(0)}=0,$ which contradicts $\|w_{(0)}\|_{W_2^1({\cal T})}=1.$ $\hfill\Box$

\medskip
Now, we are in position to prove the main result of this section.

\medskip
{\bf Theorem 6. }{\it The boundary value problem ${\cal B}$ has a unique solution $y\in W_2^1({\cal T}_\tau)\cap W_2^2(\widetilde{\cal T}).$
Moreover, there exists $C_2$ such that}
\begin{equation}\label{4.8}
\|y\|_{W_2^1({\cal T}_\tau)}\le C_2\|\varphi\|_{W_2^1[-\tau,0]}.
\end{equation}

\medskip
{\it Proof.} Consider the function $\Phi=[\Phi_1,\ldots,\Phi_m]\in W_2^1({\cal T}_\tau)$ determined by the formulae
$$
\Phi_1(t)=\left\{\begin{array}{cc}
\varphi(t), & -\tau\le t\le0,\\[2mm]
\displaystyle\frac{T_1-\tau-t}{\,T_1-\tau}\varphi(0),& 0< t\le T_1-\tau,\\[4mm]
0,& T_1-\tau< t\le T_1,
\end{array}\right.
\qquad \Phi_j(t)\equiv0, \quad j=\overline{2,m}.
$$
By virtue of Lemma~3, for a function $y\in W_2^1({\cal T}_\tau)$ obeying conditions (\ref{2.2-1})--(\ref{2.4}), to be a solution of the
problem ${\cal B},$ it is necessary and sufficient to satisfy (\ref{3.1}). In other words, $y$ is a solution of ${\cal B}$ if and only if
$x:=y-\Phi\in {\cal W}$ (which is equivalent to (\ref{2.2-1})--(\ref{2.4})) and
\begin{equation}\label{4.9}
B(\Phi,w)+B(x,w)=0 \quad \forall w\in {\cal W}
\end{equation}
(which, in turn, is equivalent to (\ref{3.1})).

Since $B(w,w)={\cal J}(w),$ Lemma~6 implies that $(\,\cdot\,,\,\cdot\,)_{\cal W}:=B(\,\cdot\,,\,\cdot\,)$ is
an inner product in~${\cal W}.$ Moreover, we have the estimate
\begin{equation}\label{4.10}
|B(\Phi,w)|=\alpha_1\Big|\int_0^{T_1}\ell_1\Phi(t)\ell_1w(t)\,dt\Big| \le
M\|\varphi\|_{W_2^1[-\tau,0]}\|w\|_{\cal W},
\end{equation}
where $\|w\|_{\cal W}=\sqrt{(w,w\,)_{\cal W}}.$ Thus, by virtue of the Riesz theorem on the general form of a linear bounded functional in a
Hilbert space, there exists a unique $x\in {\cal W}$ such that (\ref{4.9}) is fulfilled. Hence, the problem ${\cal B}$ has the unique
solution $y=\Phi+x.$

Finally, according to (\ref{4.9}) and (\ref{4.10}), we have
$$
\|x\|_{\cal W}\le M\|\varphi\|_{W_2^1[-\tau,0]},
$$
which along with Lemma~6 gives
$$
\|x\|_{W_2^1({\cal T}_\tau)}=\|x\|_{W_2^1(\widetilde{\cal T})}\le
\frac{M}{\sqrt{C_1}}\|\varphi\|_{W_2^1[-\tau,0]}.
$$
Using also the estimate
$$
\|\Phi\|_{W_2^1({\cal T}_\tau)}^2 =\|\varphi\|_{W_2^1[-\tau,0]}^2 +\|\Phi_1\|_{W_2^1[0,T_1-\tau]}^2
=\|\varphi\|_{W_2^1[-\tau,0]}^2 +\frac{(T_1-\tau)^2+3}{3(T_1-\tau)}\varphi^2(0)\le
M_1\|\varphi\|_{W_2^1[-\tau,0]}^2,
$$
we arrive at (\ref{4.8}).
% with $C_2=2M/C_1.$
$\hfill\Box$

\medskip
{\bf Remark 4.} By virtue of (\ref{4.0}) and (\ref{4.2}), the inner product $(\,\cdot\,,\,\cdot\,)_{\cal W}$
is equivalent to $(\,\cdot\,,\,\cdot\,)_{W_2^1(\widetilde{\cal T})}.$

\medskip
{\bf Remark 5.} Theorem~5 and Theorem~6 also generalize the corresponding results in \cite{Skub-94} for
(\ref{1.1}) though with $a=0$ but to the case of stepwise coefficients $b$ and $c$ when the distances between
all neighboring discontinuities including the ends of the interval is greater than $\tau.$
\\

{\bf 6. Stochastic interpretation}
\\

{\bf 6.1.} As demonstrated in Section~2, an extension of equation (\ref{1.1}) to a graph can be considered as
the situation when the coefficients in (\ref{1.1}) may randomly change their values at certain moments of
time. In order to give a systematic treatment of this view, one can employ the theory of stochastic
processes, which have applications in many fields of science and technology such as biology, chemistry,
physics, image processing, signal processing, control theory, information theory, computer science, ecology,
finance, neuroscience and telecommunications.

In general, a stochastic process can be defined as a collection of random variables $b(t;\omega)$ on a common
probability space $(\Omega, {\cal F}, P)$ and indexed by the parameter $t\in{\mathbb T}$ (see, e.g.,
\cite{Lamp}). Here, $\Omega$ is a sample space, ${\cal F}$ is a $\sigma$-algebra of events on $\Omega,$ and
$P$ is a probability measure. All random variables $b(t;\omega)$ must have a common range-space ${\cal S},$
called the process {\it state space}.

The set ${\mathbb T}$ is called {\it index set} or {\it parameter set} and usually associated with time. In
this connection, there exist {\it discrete-time} and {\it continuous-time} stochastic processes according to
whether the corresponding index set is discrete or it is an interval of the real line, respectively.

Any stochastic process can be considered as a function of two variables $t$ and $\omega.$ For each fixed
$t\in{\mathbb T},$ we have a concrete random variable $b(t;\omega),$ while fixing $\omega\in\Omega$ gives a
usual function of $t.$ The latter is called a {\it trajectory} or a {\it path-function} or a {\it
realization} of the process.

Considering various realizations of a process as values of a single random variable gives rise to an
alternative reference to stochastic processes as random functions. They also can be referred to as random
sequences or random vectors if their index sets are countable or finite, respectively.

The process realizations in the latter case have the form
$$
[b(t_1;\omega),\ldots,b(t_n;\omega)], \quad \omega\in\Omega,
$$
where the points $t_1<t_2<\ldots<t_n$ form the corresponding index set ${\mathbb T}.$

When a discrete index set ${\mathbb T}$ is naturally embedded into ${\mathbb R}$ while the set
${\mathbb T}\cap(0,T)$ is finite for any $T>0,$  all realizations of the corresponding stochastic process
can be viewed as stepwise functions on each interval $(0,T).$

The tree structure given in Section~3 allows the coefficients in (\ref{1.1}) to be various discrete-time
processes with a finite number of states in ${\mathbb R}.$ Countably many states can be covered as well but
by considering a tree with infinitely many edges emanating from each internal vertex.

\medskip
{\bf 6.2.} For an illustration, consider first the control system determined by the equation
\begin{equation}\label{6.1}
y'(t)+by(t)=u(t), \quad 0<t<T,
\end{equation}
which from the state $y(0)=\varphi\in{\mathbb R},$ is to be brought into the state $y(T)=0$ for certain $T>1.$

Assume that $u(t)$ is a real-valued function in $L_2(0,T),$ while $b$ now can take two different real values
$\theta_1$ and $\theta_2$ with probabilities $p$ and $1-p,$ respectively, so that any change between them can
occur only at the discrete times $t\in{\mathbb N}.$ For definiteness, also let $b=\theta_1$ for $t<1.$ For
example, $p=1$ would lead to a usual equation (\ref{6.1}) with the constant coefficient $b\equiv\theta_1.$
When $p\in(0,1),$ the coefficient $b$ is a Bernoulli-type process, for which we have ${\mathbb T}={\mathbb
N}.$

Under such settings, a solution $y$ of equation (\ref{6.1}) can be understood as a usual solution of a system
of differential equations on a binary tree  ${\cal T}.$

Indeed, for $t\in(0,1),$ the function $y_1(t):=y(t)$ is a solution of the usual equation
\begin{equation}\label{6.1-0}
y_1'(t)+\theta_1y_1(t)=u_1(t), \quad 0<t<1,
\end{equation}
where we also have put $u_1(t):=u(t)$ on the interval $(0,1),$ which, in turn, is assumed to parametrize the
first edge $e_1=[v_0,v_1]$ of ${\cal T}.$ For $1<t<2,$ we have two different equations
\begin{equation}\label{6.1-1}
\tilde y_2'(t)+\theta_1\tilde y_2(t)=u(t), \quad \tilde y_3'(t)+\theta_2\tilde y_3(t)=u(t), \quad 1<t<2,
\end{equation}
where $\tilde y_j(t)$ stands for $y(t)$ on $[1,2]$ when $b=\theta_{j-1}$ for $t\in[1,2).$ Clearly,
\begin{equation}\label{6.1-2}
y_1(1)=\tilde y_j(1),\quad j=2,3.
\end{equation}
Denote $y_j(t):=\tilde y_j(t+1),\;j=2,3,$ and rewrite (\ref{6.1-1}) in the form
\begin{equation}\label{6.1-3}
y_2'(t)+\theta_1 y_2(t)=u_2(t), \quad y_3'(t)+\theta_2 y_3(t)=u_3(t), \quad 0<t<1,
\end{equation}
assuming that $y_j(t)$ is defined on the edge $e_j=[v_1,v_j]$ of ${\cal T}.$ Both $u_2(t)$ and $u_3(t)$ can
stand for one and the same $u(t+1)$ for $t\in(0,1).$ However, when focusing on the control problem on a tree,
it is reasonable to allow $u_2(t)$ and $u_3(t)$ to be different functions in $L_2(0,1).$

The relations (\ref{6.1-2}) become the continuity conditions for the system (\ref{6.1-0}), (\ref{6.1-3}) at $v_1:$
$$
y_1(1)=y_2(0)=y_3(0).
$$

Further, assuming that $T\in(k-1,k]$ for some natural $k\ge2,$ and continuing this process, we arrive at a
system of differential equations on the resulting binary tree ${\cal T}$ of height $T,$ whose edges
correspond to the different outcomes of~$b:$
\begin{equation}\label{6.2}
\ell_jy(t):=y_j'(t)+b_jy_j(t)=u_j(t),\quad 0<t< T_j, \quad j=\overline{1,m},
\end{equation}
where $y_j(t)$ and $u_j(t)$ are defined on the edge $e_j=[v_{k_j},v_j].$ The other settings of Section~3 take
the form: $m=2^k-1,$ $d=2^{k-1}-1,$ $T_j=1,\,j=\overline{1,d},$ $T_j=T-k+1,\,j=\overline{d+1,m},$ and
\begin{equation}\label{6.2-1}
k_j=\frac12\left\{\begin{array}{cl}j & {\rm for \;\; even}\;\;j,\\
j-1 & {\rm for \;\;odd}\;\;j,
\end{array}\right. \quad b_j=\left\{\begin{array}{cl}\theta_1 & {\rm for \;\;odd}\;\;j,\\
\theta_2 & {\rm for \;\; even}\;\;j,
\end{array}\right. \;\; j=\overline{1,m}.
\end{equation}
In addition to (\ref{6.2}), there are also the continuity conditions (\ref{2.2-1}) at the internal vertices
of~${\cal T}.$

Obviously, all simple paths ${\cal E}_j:=\{e_{k_j^{<\nu>}}\}_{\nu=\overline{0,k-1}},$ $j=\overline{d+1,m},$
in ${\cal T},$ connecting the related boundary vertex $v_j$ with the root $v_0,$ give all realizations of $b$
in (\ref{6.1}) on the horizon~$(0,T).$

The above control problem means bringing the system (\ref{2.2-1}), (\ref{6.2}) from the initial state
\begin{equation}\label{6.3}
y_1(0)=\varphi
\end{equation}
into the state
\begin{equation}\label{6.4}
y_j(T_j)=0, \quad j=\overline{d+1,m}.
\end{equation}

The corresponding energy functional $J(y)$ in (\ref{2.5}) will contain $\ell_jy$ defined in (\ref{6.2}) along with the weights $\alpha_j$ that,
in accordance with the assignment of $\theta_1$ and $\theta_2$ to the edges as in~(\ref{6.2-1}), can be recurrently obtained by the
relations
\begin{equation}\label{6.4-1}
\alpha_1=1, \quad \alpha_j=\alpha_{k_j}\left\{\begin{array}{cc}p & {\rm for \;\; odd}\;\;j,\\
1-p & {\rm for \;\; even}\;\;j,
\end{array}\right. \quad j=\overline{2,m},
\end{equation}
(see Fig.~5). Obviously, the probability of the path ${\cal E}_j$  equals $\alpha_j,$ and
\begin{equation}\label{6.5}
\sum_{j=d+1}^m \alpha_j=1.
\end{equation}
\begin{center}
\unitlength=0.75mm
\begin{picture}(167,80)

 \put(0,45){\line(1,0){30}}
 \put(30,45){\line(1,1){21.2}}
 \put(30,45){\line(1,-1){21.2}}
 \put(51.2,66.2){\line(3,1){28.5}}      \put(120,75.8){\line(1,0){30}}   \put(120,75.8){\line(2,-1){27}}
 \put(51.2,66.2){\line(3,-1){28.5}}
 \put(51.2,23.8){\line(3,1){28.5}}      \put(120,45){\line(1,0){30}}
 \put(51.2,23.8){\line(3,-1){28.5}}     \put(120,14.2){\line(1,0){30}}

 \multiput(95,45)(2.5,0){3}{\circle*{0.8}}
 \multiput(133,29.6)(2.5,0){3}{\circle*{0.7}}

 \put(0,45){\circle*{1}}
 \put(30,45){\circle*{1}}
 \put(51.2,66.2){\circle*{1}}
 \put(51.2,23.8){\circle*{1}}
 \put(79.7,75.8){\circle*{1}}
 \put(79.7,56.6){\circle*{1}}
 \put(79.7,33.4){\circle*{1}}
 \put(79.7,14.2){\circle*{1}}

 \put(-8,44){\small $v_0$}
 \put(25,41){\small $v_1$}
 \put(47,69){\small $v_2$}
 \put(47,19.7){\small $v_3$}
 \put(80,74.5){\small $v_4$}
 \put(80,55.3){\small $v_5$}
 \put(80,32.1){\small $v_6$}
 \put(80,12.9){\small $v_7$}

 \put(150,75.8){\circle*{1}}
 \put(147,62.3){\circle*{1}}
 \put(150,45){\circle*{1}}
 \put(150,14.2){\circle*{1}}

 \put(120,75.8){\circle*{1}}
 \put(120,45){\circle*{1}}
 \put(120,14.2){\circle*{1}}

 \put(152,74.7){\small $v_{d+1}$}
 \put(149,61.2){\small $v_{d+2}$}
 \put(152,44.9){\small $v_{d+3}$}
 \put(152,12.9){\small $v_m$}

 \put(106,74.7){\small $v_{k_{d+1}}$}
 \put(106,45){\small $v_{k_{d+3}}$}
 \put(112,12.9){\small $v_d$}

 \put(8.5,47){\tiny $\alpha_1=1$}
 \put(24.5,60){\tiny $\alpha_2=1-p$}
 \put(29.5,30){\tiny $\alpha_3=p$}
 \put(69,70){\tiny $\alpha_4=(1-p)^2$}
 \put(69,62){\tiny $\alpha_5=(1-p)p$}
 \put(69,27.6){\tiny $\alpha_6=(1-p)p$}
 \put(69,19.6){\tiny $\alpha_7=p^2$}

 \put(119,78){\tiny $\alpha_{d+1}=(1-p)^{k-1}$}
 \put(138,68){\tiny $\alpha_{d+2}=(1-p)^{k-2}p$}
 \put(118,47.6){\tiny $\alpha_{d+3}=(1-p)^{k-2}p$}
 \put(124.5,16.5){\tiny $\alpha_m=p^{k-1}$}

\put(14,0){\small Fig. 5. Distribution of the probabilities on the binary tree ${\cal T}$}

\end{picture}
\end{center}

In accordance with Theorems~5 and~6, the corresponding variational problem (\ref{2.5}) under the conditions
(\ref{2.2-1}), (\ref{6.3}) and (\ref{6.4}) has a unique solution $[y_1,\ldots, y_m],$ which also uniquely
solves the boundary value problem for the equations
\begin{equation}\label{6.5-1}
y_j''(t)=b_j^2y_j(t), \quad 0< t< T_j, \quad j=\overline{1,m},
\end{equation}
under the conditions (\ref{2.2-1}), (\ref{6.3}), (\ref{6.4}) as well as the Kirchhoff conditions (\ref{3.14})
with all $\gamma_j=0.$

\medskip
{\bf Example 4.} Let $\theta_1=\theta_2.$ Then all $b_j$ are equal to $\theta_1,$ while all the paths ${\cal
E}_j,\;j=\overline{d+1,m}$ become artificial copies of one and the same realization of $b$ on the time
horizon $(0,T).$ Moreover, the optimal trajectory $[y_1,\ldots,y_m]$ becomes independent of $p.$ Indeed, let
$p=1/2.$ Then, by (\ref{6.2-1}) and (\ref{6.4-1}), the Kirchhoff conditions (\ref{3.14}) take the form
\begin{equation}\label{6.5-2}
2y_j'(1)=y_{2j}'(0)+y_{2j+1}'(0), \quad j=\overline{1,d}.
\end{equation}
Analogously to Example~2, due to the symmetry of the boundary value problem (\ref{2.2-1}), (\ref{6.3}),
(\ref{6.4}), (\ref{6.5-1}), (\ref{6.5-2}), we conclude that the optimal trajectory on ${\cal T}$ have the
groups of equal components
$$
y_{2^j}(t)\equiv y_{2^j+1}(t)\equiv\ldots\equiv y_{2^{j+1}-1}(t), \quad j=\overline{1,k-1}.
$$
In other words, the composite trajectory $[y_{k_j^{<\nu>}}]_{\nu=\overline{0,k-1}}$ is independent of
$j\in\{d+1,\ldots,m\}.$ It can be easily seen that the optimal trajectory $[y_1,\ldots,y_m]$ for any other
$p\in[0,1]$ will be the same. Hence, the problem (\ref{2.2-1}), (\ref{6.3}), (\ref{6.4}), (\ref{6.5-1}),
(\ref{6.5-2}) can be reduced to the problem on an interval
$$
y''(t)=b^2y(t), \quad 0<t<T, \quad y(0)=\varphi, \quad y(T)=0,
$$
%whose unique solvability can be checked directly.
while the connection between $y(t)$ and $y_{2^j}(t)$ for $j=\overline{0,k-1}$ can be expressed as
$$
y(t)=y_{2^j}(t-j), \;\; t\in(j,j+1], \;\; j=\overline{0,k-2}, \quad y(t)=y_{d+1}(t-k+1), \;\; t\in(k-1,T).
$$

\medskip
{\bf 6.3.} In general, one can assume all the coefficients in equation (\ref{1.1}) to be discrete-time
stochastic processes with finite numbers of states. The case of infinitely many states also can be covered.
%Different stochastic coefficients can have different as well as non-uniformly spaced discrete parameter sets,
%which always can be combined into one common index set ${\mathbb T}.$ Therefore, any such situation can be
%put into the tree structure described in Section~3.

Although Theorems~5 and~6 hold for any positive weights $\alpha_j,$ the stochastic interpretation implies
their special choice in accordance with the probabilities of the corresponding outcomes as was made, in
particular, in (\ref{1.11-1}) or (\ref{6.4-1}). In general, put $p_1=1$ and let $p_j\in(0,1)$ for
$j=\overline{2,m}$ be the conditional probability of the scenario corresponding to the edge $e_j$ provided
that the scenario of the edge $e_{k_j}$ is already known to be fulfilled. Thus, we have
$$
\sum_{\nu\in V_j}p_\nu=1, \quad j=\overline{0,d}.
$$
Then naturally generalizing the definitions (\ref{1.11-1}) and (\ref{6.4-1}), we recurrently put
\begin{equation}\label{6.6}
\alpha_1=1, \quad \alpha_j=\alpha_{k_j}p_j, \quad j=\overline{2,m}.
\end{equation}
This choice of $\alpha_j$ will be referred to as {\it probabilistic}. Obviously, it implies (\ref{6.5}). That
is, the total probability of all paths in ${\cal T}$ from the root to the remaining boundary vertices
equals~$1.$

\medskip
Finally, we show that one needs to solve the variational problem for any tree ${\cal T}$ only once. In other
words, after arriving at any internal vertex, there is no sense to look for a new optimal trajectory that
would meet the remaining subtree since it will surely coincide with the restriction to this subtree of the
optimal trajectory obtained for the entire ${\cal T}$ whenever the weights $\alpha_j,$ $j=\overline{1,m},$
are fixed from the very beginning.

This can be illustrated using Example~3, in which the parts $y_2$ and $y_3$ of the optimal trajectory
$[y_1,y_2,y_3]$ cannot be improved since they are already straight lines connecting the state taken at the
internal vertex~$v_1$ with the equilibrium positions at both possible final destinations $v_2$ and $v_3,$
respectively.

In general, let ${\cal T}^\mu$ for some $\mu\in\{2,\ldots,m\}$ be the subtree of ${\cal T}$ that has the root
$v_{k_\mu}.$ Then the control problem (\ref{2.1})--(\ref{2.4}) on ${\cal T}$ can be naturally restricted to
${\cal T}^\mu$ with the prehistory
$$
y_\mu(t)=y_{k_\mu}(t+T_{k_\mu}), \quad t\in[-\tau,0],
$$
since the edge $e_{k_\mu}$ does not belong to ${\cal T}^\mu.$ Speaking of this restriction, we assume that
$y_{k_\mu}(t)$ is the component of the optimal trajectory $y=[y_1,\ldots,y_m]$ for the original problem
(\ref{2.1})--(\ref{2.4}).

\medskip
{\bf Proposition 1. }{\it Let the weights $\alpha_j,\,j=\overline{1,m},$ be fixed and $\mu\in\{2,\ldots,m\}.$
Then the optimal trajectory for the restriction of the control problem (\ref{2.1})--(\ref{2.4}) to the
subtree ${\cal T}^\mu$ coincides with the restriction $y|_{{\cal T}^\mu}$ to ${\cal T}^\mu$ of the optimal
trajectory $y$ for this problem on ${\cal T}.$ }

\medskip
{\it Proof.} Let ${\cal W}_\mu$ and $B_\mu(y,v)$ have the same sense for ${\cal T}^\mu$ as ${\cal W}$ and
$B(y,v),$ respectively, have for ${\cal T}$ (see Lemma~1). Then the space ${\cal W}_\mu$ is a subspace of
${\cal W}.$ Specifically, ${\cal W}_\mu$ consists of the restrictions $v^0|_{{\cal T}^\mu}$ to ${\cal
T}^\mu_\tau$ of all the functions $v^0$ in ${\cal W}$ that vanish on ${\cal T}\setminus{\cal T}^\mu.$

Then by Lemma~1, we have
$$
B_\mu(y|_{{\cal T}^\mu},v^0|_{{\cal T}^\mu})=B(y,v^0)=0 \quad \forall v^0\in{\cal W}.
$$
Thus, it remains to use the analog of Lemma~1 for ${\cal T}^\mu.$ $\hfill\Box$

\medskip
{\bf Corollary 1. }{\it If the conditional probabilities $p_j,\;j=\overline{2,m}$ are fixed, then the
assertion of Proposition~1 remains true under the probabilistic choices of the weights both for ${\cal T}$
and for ${\cal T}^\mu.$}

\medskip
{\it Proof.} Let ${\cal J}_\mu(z)$ be the functional for the restriction of the control problem
(\ref{2.1})--(\ref{2.4}) to the subtree ${\cal T}^\mu$ involving those weights $\alpha_j$ that are
probabilistic
% for the functional ${\cal J}(y)$ defined by (\ref{2.5})
for ${\cal T}.$ Denote by $\widetilde{\cal J}_\mu(z)$ the functional for~${\cal T}^\mu$ with its own
probabilistic weights. Then, in accordance with (\ref{6.6}), we have ${\cal
J}_\mu(z)=\alpha_{k_\mu}p_\mu\widetilde{\cal J}_\mu(z).$ Hence, ${\cal J}_\mu(z)$ and $\widetilde{\cal
J}_\mu(z)$ can be minimized only simultaneously. $\hfill\Box$

\bigskip
%{\bf Acknowledgements.} The author is grateful to Professor A.\,L. Skubachevskii for his suggestion to
%consider the problem of damping a control system with aftereffect in connection with the idea of global delay
%proposed by the author, and to Dr. Maria Kuznetsova for reading the manuscript and making useful comments, as
%well as the anonymous referees for their valuable recommendations.
%
%\medskip
{\bf Funding.} This research was supported by Russian Science Foundation, Grant No. 22-21-00509, https://rscf.ru/project/22-21-00509/

The author's affiliations:

\begin{itemize}
\item{} Saratov State University, Saratov, 410012 Russia

\item Moscow Center for Fundamental and Applied Mathematics, Moscow, 119991 Russia

\item Lomonosov Moscow State University, Moscow, 119991 Russia
\end{itemize}

{\it E-mail address}: {\color{blue}\bf buterinsa@sgu.ru}

\end{document}